\documentclass[12pt,reqno]{amsart}
\usepackage{amssymb,amsmath}
\newtheorem*{thm}{Theorem}
\newtheorem*{cla}{Claim}

\newcommand{\diff}{\mathop{\rm Diff\mspace{2mu}}}
\newcommand{\id}{\mathop{\rm Id}}
\newcommand{\N}{{\mathbb N}}
\newcommand{\Z}{{\mathbb Z}}
\newcommand{\R}{{\mathbb R}}
\newcommand{\al}{{\alpha}}
\newcommand{\be}{{\beta}}
\newcommand{\Al}{{\mathcal A}}

\renewcommand{\phi}{\varphi}
\newcommand{\ree}{\mathop{\rm\,Re\,}}

\newcommand{\rrr}[1]{{\rm (\ref{#1})}}

\title
{Slow area-preserving diffeomorphisms of the torus}
\author
{Alexander Borichev}
\keywords{Area-preserving diffeomorphisms, growth sequences.\newline
$2000$ {\it Mathematical Subject Classification} 37C05, 26A18, 28D05, 58D05.
\newline{\textit{\small E-mail}: \texttt{\small borichev@math.u-bordeaux.fr}}}

\begin{document}
\begin{abstract} We construct area-preserving 
real analytic diffeomorphisms of the torus with 
unbounded growth sequences of arbitrarily slow growth.
\end{abstract}
\maketitle

Given a smooth compact manifold $M$, consider the group
$\diff(M)$ of diffeomorphisms of $M$. For every $f\in \diff(M)$
we define the growth sequence of $f$:
$$
\Gamma_n(f)=\max\bigl(\max_{x\in M}\|d_xf^n\|,\max_{x\in M}\|d_xf^{-n}\|\bigr),
\qquad n\in\N,
$$
where $f^n$ is the $n$-th iteration of $f$, $f^{-n}$ is the 
$n$-th iteration of $f^{-1}$, and $\|d_xf\|$ is the operator
norm of the differential of $f$ at the point $x\in M$. 
Conjugations of $f$ in the group
$\diff(M)$ generate equivalent growth sequences:
$$
c(g)\Gamma_n(g^{-1}fg)\le \Gamma_n(f)\le C(g)\Gamma_n(g^{-1}fg), \qquad 
g\in \diff(M),\,
n\in\N.
$$

The asymptotics of the growth sequence is a basic dynamic invariant 
(see \cite{KH}). D'Ambra and Gromov \cite[7.10.C]{AG} proposed to study the behavior
of growth sequences for various classes of diffeomorphisms. In particular, 
it is interesting to find examples of unbounded growth sequences of
slow growth (see also the references in \cite[7.10.C]{AG}). 
We call the diffeomorphisms generating such growth sequences
the slow diffeomorphisms. 

Recently, Polterovich and Sodin \cite{PS} obtained several results on
the growth sequences of smooth order-preserving diffeomorphisms of the interval
$[0,1]$. In particular, they proved \cite[Theorem~1.7]{PS} that for any sequence 
$\{a_n\}$ of positive numbers tending to infinity, there exists a $C^\infty$-smooth
diffeomorphism $f$, $f\ne\id$, 
such that 
$$
\liminf_{n\to\infty}\frac{\Gamma_n(f)}{a_n}\le 1.
$$
On the other hand, a simple argument (see \cite{PS}) shows that
for any diffeomorphism $f\ne\id$,
$$
\sum_{n\ge 1}\frac{1}{\Gamma_n(f)}<\infty.
$$

Furthermore, Polterovich proved \cite[Theorem 1.3]{P} that
for every $0<\beta<1$, 
there exists an area-preserving real analytic 
diffeomorphism $f$ of the torus such that
\begin{gather*}
\Gamma_n(f)\le Cn^{\beta}\log n,\\
\limsup_{n\to\infty}\frac{\Gamma_n(f)}{n^\beta}>0.
\end{gather*}

In this note, we improve somewhat the result by  Polterovich
by producing area-preserving 
real analytic diffeomorphisms of the torus with arbitrarily slowly growing
unbounded growth sequences.

\begin{thm} Let $\phi$ be a positive increasing {\rm(}unbounded{\rm)} function 
on $\R_+$ such that
$\phi(x)=o(x)$, $x\to\infty$. There exists an area-preserving 
real analytic diffeomorphism $f$ of the torus such that
\begin{equation}
\left.
\begin{gathered}
\Gamma_n(f)\le \phi(n),\\
\limsup_{n\to\infty}\frac{\Gamma_n(f)}{\phi(n)}>0.
\end{gathered}
\right\}
\label{1}
\end{equation}
\end{thm}

For some related questions on the asymptotics of diffeomorphisms
with fixed points see \cite{P}. Other recent results on the behavior
of the growth sequences are in \cite{B}, \cite{PS}.

\begin{proof}[Proof of the Theorem] 
We represent the torus as the product $[0,1)\times[0,1)$,
and define, as in \cite{P},
$$
f(x,y)=(\{x+\al\},\{y+cF(x)\}),\qquad x,y\in[0,1),
$$
for $\al\in\R$, $c\in(0,+\infty)$, and a real analytic 
$1$-periodic function $F:\R\to\R$; here $\{\cdot\}$ stands
for the fractional part. Then
$$
d_xf=
\left(
\begin{matrix}
1&0\\c\,F'(x)&1
\end{matrix}
\right),
$$
and $f$ is an area-preserving 
real analytic diffeomorphism of the torus.

Define the Weyl sums
\begin{equation}
W(N,x,\al)=\sum_{n=0}^{N-1}F'(x+n\al).
\label{1.5}
\end{equation}
We have
$$
d_xf^N=
\left(
\begin{matrix}
1&0\\c\,W(N,x,\al)
&1
\end{matrix}
\right),\qquad N\ge 1,
$$
and 
$$
d_xf^{-N}=
\left(
\begin{matrix}
1&0\\-c\,W(N,x-N\al,\al)&1
\end{matrix}
\right),\qquad N\ge 1.
$$

Therefore, for \rrr{1} to hold it suffices that $F$ and $\al$
satisfy the following condition:
\begin{equation}
0<\limsup_{N\to\infty}\max_{0\le x<1}\frac{W(N,x,\al)}{\phi(N)}<\infty.
\label{2}
\end{equation}

Up to now our proof repeated that of Polterovich in \cite{P}.
The main difference of our argument is in the way of estimating the
Weyl sums \rrr{1.5}.

We are going to choose a sequence $\{q_k\}_{k\ge 1}$, $q_1=1$,
\begin{equation}
\frac{q_{k+1}}{100 q_k}\in\N,\qquad k\ge 1,
\label{3}
\end{equation}
and a sequence $\{r_k\}_{k\ge 1}$,
\begin{equation}
0<r_k<\exp(-q_k),\qquad k\ge 1, 
\label{4}
\end{equation}
and define
$$
F(x)=\sum_{k\ge 1}\frac{r_k}{2\pi q_k}\sin[2\pi q_k x],\qquad x\in\R\,.
$$
Then $F$ is real analytic  and $1$-periodic,
$$
F'(x)=\sum_{k\ge 1}r_k\cos[2\pi q_k x],\qquad x\in\R\,.
$$

For $\alpha\in\R$ denote
\begin{align}
\Delta_k(N,\al)&=\sum_{n=0}^{N-1}e^{2\pi iq_k n\al},\label{9.5}\\
S(N,\al)&=\sum_{k\ge 1}r_k\ree \Delta_k(N,\al),\notag\\
T(N,\al)&=\sum_{k\ge 1}r_k|\Delta_k(N,\al)|.\notag
\end{align}
Then
\begin{align*}
\Delta_k(N,\al)&=\frac{1-e^{2\pi iq_k N\al}}{1-e^{2\pi iq_k \al}},
\qquad q_k \al\in\R\setminus\Z,
\\
W(N,x,\al)&=\sum_{k\ge 1}r_k\ree \Bigl[e^{2\pi iq_k x}
\sum_{n=0}^{N-1}e^{2\pi iq_k n\al}\Bigr]\\
&=\sum_{k\ge 1}r_k\ree \Bigl[e^{2\pi iq_k x}\Delta_k(N,\al)\Bigr],
\end{align*}
and the property \rrr{2} follows from the inequalities
\begin{align}
\limsup_{N\to\infty}\frac{S(N,\al)}{\phi(N)}&>0,\label{6}\\
\limsup_{N\to\infty}\frac{T(N,\al)}{\phi(N)}&<\infty.\label{7}
\end{align}

To get \rrr{6} and \rrr{7}, we should first study the behavior 
of the sums $\Delta_k(N,\alpha)$. Essentially, if 
the fractional part of $q_k\alpha$ is of order $1/M$, then 
$\Delta_k(N,\alpha)$ behaves like $N/M$ for $N$ smaller than $M$, 
and is bounded by a constant times $M$ for all $N$. After that, 
in an inductive process we approximate $\phi$ from below on 
an infinite sequence of points by a weighted sum of 
$\Delta_k(N,\alpha)$ with a lacunary sequence $q_k$ and
a suitable $\alpha$.

Our first observation is as follows. Fix $n\ge 1$, 
suppose that the numbers $q_1,\ldots,q_{n+1}$ satisfy condition (\ref{3}),
and define 
$$
k_n=\sum_{1\le s\le n}\frac{q_n}{q_s}\in\N.
$$
 
\begin{cla} Suppose that $\be$ belongs to the interval
\begin{equation}
\Al_n=\Bigl\{\be: \frac{q_n}{q_{n+1}} \le q_n\be-k_n\le \frac{2q_n}{q_{n+1}} \Bigr\}.
\label{nn}
\end{equation}
Then
\begin{align}
\Bigl|\frac{\Delta_n(N,\be)}{N}-1\Bigr|&\le \frac 12,\qquad 1\le N\le 
\frac{q_{n+1}}{100q_n},\label{8}\\
|\Delta_n(N,\be)|&\le \frac{q_{n+1}}{q_n},\qquad N\ge 1.
\label{9} 
\end{align}
\end{cla}

\begin{proof}
Applying the Taylor formula to the function $x\mapsto \exp ix$,
and using that $\be\in\Al_n$ and
$q_n/q_{n+1}\le 1/100$, we obtain 
\begin{equation}
\biggl|\frac{e^{2\pi i (q_n\be-k_n)}-1}{2\pi i (q_n\be-k_n)}-1\biggr|\le
\biggl|\frac{(2\pi)^2 (q_n\be-k_n)^2}
{4\pi (q_n\be-k_n)}\biggr|\le \frac1{10}.\label{10}
\end{equation}
Furthermore, if $1\le N\le q_{n+1}/(100q_n)$, then
\begin{multline*}
\biggl|\frac{e^{2\pi i N(q_n\be-k_n)}-1}{2\pi i N(q_n\be-k_n)}-1\biggr|\le
\biggl|\frac{(2\pi)^2 N^2(q_n\be-k_n)^2}
{4\pi N(q_n\be-k_n)}\biggr|\\
=\pi N(q_n\be-k_n)\le\frac1{10}.
\end{multline*}
Hence,
$$
\biggl|\frac{\Delta_n(N,\be)}{N}-1\biggr|=
\biggl|\frac{1-e^{2\pi i N(q_n\be-k_n)}}
{N(1-e^{2\pi i (q_n\be-k_n)})}-1\biggr|\le
\frac 12,\qquad 1\le N\le \frac{q_{n+1}}{100q_n},
$$
and \rrr{8} is proved. 

Next,
$$
|1-e^{2\pi i Nq_n\be}|\le 2,\qquad N\ge 1.
$$
Therefore,
$$
|\Delta_n(N,\be)|=
\biggl|\frac{1-e^{2\pi i Nq_n\be}}
{1-e^{2\pi i q_n\be}}\biggr|\le
\frac 2{|1-e^{2\pi i (q_n\be-k_n)}|}.
$$
Using \rrr{10}, we get
$$
|\Delta_n(N,\be)|\le
\frac 3{2\pi(q_n\be-k_n)}\le 
\frac{q_{n+1}}{q_n},\qquad N\ge 1,
$$
and \rrr{9} is proved. 
\end{proof}

Now, to obtain \rrr{6} and \rrr{7}, we define $q_k$, $r_k$, and $\al$  
in an inductive process.
Without loss of generality we assume that $\phi(1)\ge 2$.
Set $S^0(N,\be)=T^0(N,\be)=0$, $q_1=1$, $N_0=M_0=1$, $\Al_0=(0,1)$. 
On the induction step $p\ge 1$ we start with
sequences $\{q_j\}_{1\le j\le p}$,
$\{r_j\}_{1\le j<p}$, $\{N_j\}_{0\le j<p}$, $\{M_j\}_{0\le j<p}$,
and the interval $\Al_{p-1}$ (defined by $\{q_j\}_{1\le j\le p}$),
such that for every $\be\in\Al_{p-1}$, the function $S^{p-1}$,
$$
S^{p-1}(N,\be)=\sum_{1\le n<p}r_n\ree \Delta_n(N,\be),
$$
satisfies
\begin{equation}
S^{p-1}(N_j,\be)\ge \frac 1{100}\phi(N_j)+2^{-p},\qquad 1\le j<p,
\label{11}
\end{equation}
and for every $\be\in\Al_{p-1}$ the function $T^{p-1}$,
$$
T^{p-1}(N,\be)=\sum_{1\le n<p}r_n|\Delta_n(N,\be)|,
$$
satisfies
\begin{align}
T^{p-1}(N,\be)&\le 200\,\phi(N)-2^{-p},\qquad N\ge 1,
\label{12}\\
T^{p-1}(N,\be)&\le \frac1{100}\phi(M_{p-1}),\qquad N\ge 1.
\label{13}
\end{align}

By \rrr{9.5}, for any $\be$, $N$,
\begin{equation}
|\Delta_p(N,\be)|\le N,
\label{111}
\end{equation}
and we can choose $r_p$
satisfying \rrr{4} such that for all $\be\in\Al_{p-1}$,
\begin{equation}
r_p(N+|\Delta_p(N,\be)|)\le 2^{-p-1},\qquad 1\le N\le \max(M_{p-1},N_{p-1}).
\label{14}
\end{equation}

Since $\phi$ is increasing and $\lim_{N\to\infty}\phi(N)/N=0$, we can find
$$
N_p>\max(M_{p-1},N_{p-1})
$$ 
such that
\begin{gather}
r_pN\le \phi(N),\qquad N< N_p,\label{15}\\
\phi(N_p)\le r_pN_p\le \phi(N_p)+1.\label{16}
\end{gather}

Set $q_{p+1}=100q_pN_p$. 
Then we define $\Al_p\subset\Al_{p-1}$ by the formula \rrr{nn}.
In the estimates to follow
we assume that $\be\in\Al_p$, and hence, by the Claim, \rrr{8}
and \rrr{9} hold with $n=p$.

By \rrr{8} and \rrr{16},
\begin{equation}
r_p \ree\Delta_p(N_p,\be)\ge \frac12 r_pN_p
\ge \frac12 \phi(N_p).\label{17}
\end{equation}
By \rrr{111} and \rrr{15},
$$
r_p|\Delta_p(N,\be)|\le r_p N\le \phi(N),\qquad N<N_p,
$$
and by \rrr{9} and \rrr{16},
\begin{equation}
r_p|\Delta_p(N,\be)|\le 100\, r_p N_p\le
100(\phi(N_p)+1),\qquad N\ge N_p.\label{18}
\end{equation}
Thus, 
\begin{equation}
r_p|\Delta_p(N,\be)|\le 100(\phi(N)+1),\qquad N\ge 1.\label{19}
\end{equation}

Now, by \rrr{11} and \rrr{14},
$$
S^p(N_j,\be)\ge \frac 1{100}\phi(N_j)+2^{-p-1},\qquad 1\le j< p,
$$
and by \rrr{13} and \rrr{17},
$$
S^p(N_p,\be)\ge \frac 1{100}\phi(N_p)+2^{-p-1}.
$$
Thus, 
$$
S^p(N_j,\be)\ge \frac 1{100}\phi(N_j)+2^{-p-1},\qquad 1\le j\le p.
$$

Furthermore, by \rrr{12} and \rrr{14},
$$
T^p(N,\be)\le 200\,\phi(N)-2^{-p-1},\qquad 1\le N\le M_{p-1},
$$
and by \rrr{13} and \rrr{19},
$$
T^p(N,\be)\le 200\,\phi(N)-1,\qquad N>M_{p-1}.
$$
Thus, 
$$
T^p(N,\be)\le 200\,\phi(N)-2^{-p-1},\qquad N\ge 1.
$$

Finally, by \rrr{13}, \rrr{18},
and by the condition that $\phi$ is unbounded, there exists
$M_p$ such that
$$
T^p(N,\be)\le \frac 1{100}\phi(M_p),\qquad N\ge 1.
$$
This completes the induction step.

The intervals $\Al_n$ constitute a nested family,
\begin{gather*}
\bigcap_{n\ge 0} \Al_n=\{\al\},\\
\al=\sum_{k\ge 1}\frac1{q_k},
\end{gather*}
and all the inequalities in the induction process
are valid with $\be=\al$.

We have 
$$
S(N,\al)=\lim_{p\to\infty}S^p(N,\al),\qquad T(N,\al)=\lim_{p\to\infty}T^p(N,\al).
$$
Then the properties \rrr{11} and \rrr{12} imply that
\begin{align*}
\frac{S(N_j,\al)}{\phi(N_j)}&\ge\frac1{100},\qquad j\ge 1,\\
\frac{T(N,\al)}{\phi(N)}&\le 200,\qquad N\ge 1,
\end{align*}
and \rrr{6} and \rrr{7} follow.
The theorem is proved.
\end{proof}
\bigskip

The author is thankful to Leonid Polterovich and 
Misha Sodin for helpful remarks.

\bigskip

\noindent \textsf{\small Alexander Borichev, Department of Mathematics,} 
\newline
\noindent\textsf{\small University
of Bordeaux I, 351, cours de la Lib\'eration, 33405 Talence, France}

\end{document}